\documentclass[12pt,a4paper]{article}
\usepackage{mathrsfs}
\usepackage{amssymb}

\usepackage{amsmath}

\newtheorem{theorem}{Theorem}[section]
\newtheorem{definition}[theorem]{Definition}
\newtheorem{lemma}[theorem]{Lemma}

\begin{document}

\title{Gr\"{o}bner-Shirshov bases for Schreier extensions of
groups\footnote {The research is supported by the National Natural
Science Foundation of China (Grant No.10771077) and the Natural
Science Foundation of Guangdong Province (Grant No. 06025062).} }
\author{
Yuqun Chen  \\
%EndAName
\\
{\small \ School of Mathematical Sciences}\\
{\small \ South China Normal University}\\
{\small \ Guangzhou 510631}\\
{\small \ P. R. China}\\
{\small \ yqchen@scnu.edu.cn} }

\date{}

\maketitle \noindent{\bf Abstract.} In this paper, by using the
Gr\"{o}bner-Shirshov bases, we give characterizations of the
Schreier extensions of groups when the group is presented by
generators and relations. An algorithm to find the conditions of a
group to be a Schreier extension is obtained. By introducing a
special total order, we obtain the structure of the Schreier
extension by an HNN group.
\\
%\end{center}

\section{Preliminaries}

The concept of the extensions of groups was introduced by Otto
Schreier in \cite{s1} and \cite{s2}, also see  M. Hall's book
\cite{h}.
\begin{definition} \ (\cite{h},\cite{m}) \
Let $A,B,G$ be groups. By an exact sequence
$$
1\rightarrow A\rightarrow G\rightarrow B\rightarrow 1
$$
we mean that $A$ is a normal subgroup of $G$ such that $G/A\cong B$.
If this is the case, $G$ is called a Schreier extension of $A$ by
$B$.
\end{definition}

Some characterizations of Schreier extensions are included in
\cite{h} (see \cite{h}, pp218-228, specially Chapter 15.4. Defining
Relations and Extensions). It is told (\cite{h}, p228) that ``It is
difficult to determine the identities (in $A$) leading to conditions
for an extension", where the group $B$ is presented by generators
and relations.

In this paper, by using the Gr\"{o}bner-Shirshov bases with a
monomial total order, we give characterizations of the Schreier
extensions of groups when the group is presented by generators and
relations. An algorithm to find the extension conditions of a group
to be a Schreier extension is obtained. Then, we find a solution to
the above problem discussed by M. Hall in the book \cite{h}.

A Gr\"{o}bner-Shirshov bases for HNN-extension of group is given in
\cite{cz} by using the Generalized Composition-Diamond Lemma (see
Lemma \ref{l3.1}) with a special total order (not monomial). Because
this order is not monomial, we have to find a special order such
that we can use the Generalized Composition-Diamond Lemma to give
the extension conditions by an HNN group.

\ \

Now, we cite some concepts and results from the literature.

Let $k$ be a field, $k\langle X\rangle$ the free associative algebra
over $k$ generated by $X$ and $ X^{*}$ the free monoid generated by
$X$, where the empty word is the identity which is denoted by 1. For
a word $w\in X^*$, we denote the length of $w$ by $deg(w)$. Let
$X^*$ be a total ordered set. Let $f\in k\langle X\rangle$ with the
leading word $\bar{f}$. We say that $f$ is monic if $\bar{f}$ has
coefficient 1.

\begin{definition} (\cite{Sh}, see also \cite{b72}, \cite{b76}) \
Let $f$ and $g$ be two monic polynomials in \textmd{k}$\langle
X\rangle$ and $<$ a total order on $X^*$. Then, there are two kinds
of compositions:

$(1)$ If \ $w$ is a word such that $w=\bar{f}b=a\bar{g}$ for some
$a,b\in X^*$ with deg$(\bar{f})+$deg$(\bar{g})>$deg$(w)$, then the
polynomial
 $(f,g)_w=fb-ag$ is called the intersection composition of $f$ and
$g$ with respect to $w$.

$(2)$ If  $w=\bar{f}=a\bar{g}b$ for some $a,b\in X^*$, then the
polynomial $(f,g)_w=f - agb$ is called the inclusion composition of
$f$ and $g$ with respect to $w$.

\end{definition}

\begin{definition}(\cite{b72}, \cite{b76}, cf. \cite{Sh})
Let $S\subseteq$ $\textmd{k}\langle X\rangle$ with each $s\in S$
monic. Then the composition $(f,g)_w$ is called trivial modulo $S$
if $(f,g)_w=\sum\alpha_i a_i s_i b_i$, where each $\alpha_i\in k$,
$a_i,b_i\in X^{*}$ and $\overline{a_i s_i b_i}<w$. If this is the
case, then we write
$$
(f,g)_w\equiv0\quad mod(S,w)
$$
In general, for $p,q\in k\langle X\rangle$, we write
$$
p\equiv q\quad mod(S,w)
$$
which means that $p-q=\sum\alpha_i a_i s_i b_i $, where each
$\alpha_i\in k,a_i,b_i\in X^{*}$ and $\overline{a_i s_i b_i}<w$.
\end{definition}

\begin{definition} (\cite{b72}, \cite{b76}, cf. \cite{Sh}) \
We call the set $S$ with respect to the total order $<$ a
Gr\"{o}bner-Shirshov set (basis) in $k\langle X\rangle$ if any
composition of polynomials in $S$ is trivial relative to $S$.
\end{definition}

If a subset $S$ of $k\langle X\rangle$ is not a Gr\"{o}bner-Shirshov
basis, then we can add to $S$ all nontrivial compositions of
polynomials of $S$, and by continuing this process (maybe
infinitely) many times, we eventually obtain a Gr\"{o}bner-Shirshov
basis $S^{comp}$. Such a process is called the Shirshov algorithm.

A total order $>$ on $X^*$ is monomial if it is compatible with the
multiplication of words, that is, for $u, v\in X^*$, we have
$$
u > v \Rightarrow w_{1}uw_{2} > w_{1}vw_{2},  \ for \  all \
 w_{1}, \ w_{2}\in  X^*.
$$
A standard example of monomial order on $X^*$ is the deg-lex order
to compare two words first by degree and then lexicographically,
where $X$ is a linearly ordered set.

The following lemma was proved by Shirshov \cite{Sh} for the free
Lie algebras (with deg-lex ordering) in 1962 (see also Bokut
\cite{b72}). In 1976, Bokut \cite{b76} specialized the approach of
Shirshov to associative algebras (see also Bergman \cite{b}). For
commutative polynomials, this lemma is known as the Buchberger's
Theorem in \cite{bu65} and \cite{bu70}.

\begin{lemma}\label{l1}
(Composition-Diamond Lemma) \ Let $k$ be a field, $A=k \langle
X|S\rangle=k\langle X\rangle/Id(S)$ and $<$ a monomial order on
$X^*$, where $Id(S)$ is the ideal of $k \langle X\rangle$ generated
by $S$. Then the following statements are equivalent:
\begin{enumerate}
\item[(1)] $S $ is a Gr\"{o}bner-Shirshov basis.
\item[(2)] $f\in Id(S)\Rightarrow \bar{f}=a\bar{s}b$
for some $s\in S$ and $a,b\in  X^*$.
\item[(3)] $Irr(S) = \{ u \in X^* |  u \neq a\bar{s}b ,s\in S,a ,b \in X^*\}$
is a basis of the algebra $A=k\langle X | S \rangle$.
\end{enumerate}
\end{lemma}

For convenience, we identify a relation $v=h_v$ of an algebra
presented by generators and relations with the polynomial $v-h_v$ in
the corresponding free algebra.

\section{Gr\"{o}bner-Shirshov bases for Schreier extensions of
groups}

In this section, by using Gr\"{o}bner-Shirshov bases, we give the
structure of Schreier extensions of groups and an algorithm is given
to find all Schreier extensions of the groups $A$ by $B$.

Let $A,B$ be groups. By a factor set of $B$ in $A$, we mean a subset
of $A$ which is related to $B$.

Let the group $B$ be presented as semigroup by generators and
relations:
$$
B=sgp\langle Y|R \rangle
$$
where $R$ is a Gr\"{o}bner-Shirshov bases for $B$ with the deg-lex
order $<_B$ on $Y^*$.

For convenience, we can assume that $R$ is a minimal
Gr\"{o}bner-Shirshov bases in a sense that the leading monomials are
not contained each other as subwords, in particular, they are
pairwise different.  Let $R=\{v=h_v|v\in\Omega\}$, where $v$ is the
leading term of the polynomial $f_v=v-h_v$ in $k\langle Y \rangle$.
Let
$$
G=E(A,Y,a^y, (v))=sgp\langle A_1\cup Y| S\rangle
$$
where $A_1=A\backslash\{1\}, \ S=\{aa'=[aa'], v=h_v\cdot(v), \
ay=ya^y |v\in\Omega, \  a,a'\in A_1, \ y\in Y\}$,
$\{(v)|v\in\Omega\}\subseteq A$ a factor set of $B$ in $A$, $\psi_y:
\ A\rightarrow A, \ a\mapsto a^y$ an automorphism.

Let $A_1$ be a total ordered set. Then, by the deg-lex order, we
have the monomial total order $<_A$ on $A_1^*$.

We order the set $(A_1\cup Y)^*$. For each element $d$ in $(A_1\cup
Y)^*$, $d$ has a unique expression: $d=e_1y_1\cdots e_ty_te_{t+1}$,
where each $e_l\in A_1^*, \ y_l\in Y$. Suppose that
$d'=e_1'y_1'\cdots e_r'y_r'e_{r+1}'\in (A_1\cup Y)^*$. Let
\begin{eqnarray*}
&&wt(d)=(t,y_1,\cdots, y_t,e_1,\cdots,e_t,e_{t+1})\\
&&wt(d')=(r,y_1',\cdots, y_r',e_1',\cdots,e_r',e_{r+1}')
\end{eqnarray*}
We define $d\prec d'$ if $wt(d)<wt(d')$ lexicographically, using the
order of natural numbers and the following orders:
\begin{enumerate}
\item[(a)]\ $y_l<y_l'$ by the order in $Y$;
\item[(b)]\ $e_l<e_l'$ by the order $<_A$.
\end{enumerate}

It is clear that $\prec$ is a monomial total order on $(A_1\cup
Y)^*$ which extends the orders $<_A$ and $<_B$. We call such an
order the tower order.

In $S$, the possible compositions are related to the following
ambiguities:
$$
w_1=aa'a'', \ \ w, \ \ av, \ \ aa'y
$$
where $v\in\Omega, a,a',a''\in A_1, \ y\in Y$ and $w$ is an
ambiguity appeared in $R$.

It is clear that for $w_1=aa'a''$, the composition is trivial modulo
$S$ since $[[aa']a'']=[a[a'a'']]$.

For $w_1=aa'y$, since
$$
(aa'-[aa'],a'y-ya'^y)_{w_{_1}}=-[aa']y+aya'^y\\
\equiv -y[aa']^y+ya^ya'^y\equiv0 \ \ mod(S,w_{_1})
$$
the composition is trivial modulo $S$.

 For any $a\in A, \ w_1=av, \ v=y_1\cdots
y_r, \ h_v=y_1'\cdots y_t'$, we have
\begin{eqnarray*}
(ay_1-y_1a^{y_1},v-h_v\cdot(v))_{w_{_1}}&\equiv&
-y_1a^{\psi_{y_1}}y_2\cdots y_r+ah_v\cdot(v)\\
&\equiv& -y_1y_2\cdots y_ra^{\psi_{y_1}\psi_{y_2}\cdots \psi_{y_r}}+
y_1'\cdots y_t'a^{\psi_{y_1'}\cdots \psi_{y_t'}}\cdot(v)\\
&\equiv&-h_v\cdot(v) a^{\psi_{y_1}\psi_{y_2}\cdots \psi_{y_r}}+
h_va^{\psi_{y_1'}\cdots \psi_{y_t'}}\cdot(v)\\
&\equiv&-h_v((v) a^{\psi_{y_1}\psi_{y_2}\cdots \psi_{y_r}}-
a^{\psi_{y_1'}\cdots \psi_{y_t'}}\cdot(v))\\
&\equiv&h_v(a^{h_v}(v)-(v) a^{v}) \ \ mod(S,w_{_1})
\end{eqnarray*}
where $a^{\psi_{y_1}\cdots \psi_{y_r}}=a^v$ if $v=y_1\cdots y_r$.

Then, for $w_1=av, \ v\in \Omega$, the composition is trivial modulo
$S$ if and only if
\begin{eqnarray}\label{e6}
(v) a^{v}= a^{h_v}(v)
\end{eqnarray}
holds in $\bar{A}$, where $\bar{A}$ is the natural homomorphic image
of $A$ in $G$.

Since $R$ is a minimal Gr\"{o}bner-Shirshov bases, all compositions
in $R$ are only intersection ones. Now, for $w_1=w=v_1c=dv_2, \
v_1,v_2\in \Omega, \ c,d\in Y^*, \ deg(v_1) + deg(v_2) > deg (w)$,
we have,
$$
f_{v_1}c - df_{v_2} =dh_{v_2} - h_{v_1}c\equiv 0  \ \ mod(R,w)
$$
It means that there exists a $z\in Y^*$ such that
\begin{eqnarray*}
&&dh_{v_2} \equiv z \ \ mod(R,w)\\
&&h_{v_1}c \equiv z \ \ mod(R,w)
\end{eqnarray*}

Now, by invoking that
\begin{eqnarray*}
&&(v_1-h_{v_{_1}}\cdot(v_1),v_2-h_{v_{_2}}\cdot(v_2))_{w_{_1}}=
-(h_{v_{_1}}\cdot(v_1))c+d(h_{v_{_2}}\cdot(v_2))\\
&\equiv& (dh_{v_{_2}})(v_2)-(h_{v_{_1}}c)(v_1)^c \ \ \ mod(S,w_{_1})
\end{eqnarray*}
there exist $\xi_{(v_1,v_2)_{_w}}(v), \ \zeta_{(v_1,v_2)_{_w}}(v)\in
A$ such that
\begin{eqnarray}\label{e2}
(v_1-h_{v_{_1}}\cdot(v_1),v_2-h_{v_{_2}}\cdot(v_2))_{w_{_1}}\equiv
z(\xi_{(v_1,v_2)_{_w}}(v)-\zeta_{(v_1,v_2)_{_w}}(v)) \ \
mod(S,w_{_1})
\end{eqnarray}
where $\xi_{(v_1,v_2)_{_w}}(v)$ and $\zeta_{(v_1,v_2)_{_w}}(v)$ are
functions of $\{(v)|v\in\Omega\}$.

In fact by the previous formulas we have an algorithm to find the
functions $\xi_{(v_1,v_2)_{_w}}(v)$ and $
\zeta_{(v_1,v_2)_{_w}}(v)$.

Then, for $w_1=w$, the composition is  trivial modulo $S$ if and
only if
\begin{eqnarray}\label{e3}
\xi_{(v_1,v_2)_{_w}}(v)=\zeta_{(v_1,v_2)_{_w}}(v)
\end{eqnarray}
holds in $\bar{A}$.

It is reminded that, if  $z=y^{-1}\in Y, \ y\in Y$ and $R$ contains
the relation $f_v=y^{\epsilon}y^{-\epsilon}-1, \ \epsilon=\pm 1$,
then in the factor set $\{(v)|v\in\Omega\}$, $(v)=1$  in $\bar{A}$.

Therefore, we have the following theorem.

\begin{theorem}\label{sgt1}
Let  $A,B$ be groups, $B=sgp\langle Y|R \rangle$, where
$R=\{v=h_v|v\in\Omega\}$ is a minimal Gr\"{o}bner-Shirshov bases for
$B$ and $v$ the leading term of the polynomial $f_v=v-h_v$ in
$k\langle Y \rangle$. Let
$$
G=E(A,Y,a^y,(v))=sgp\langle A_1\cup Y| S\rangle
$$
where $A_1=A\backslash\{1\}, \ S=\{aa'=[aa'], \ v=h_v\cdot(v), \
ay=ya^y |v\in\Omega, \ a,a'\in A_1, \ y\in Y\}$, $\psi_y: \
A\rightarrow A, \ a\mapsto a^y$ an automorphism,
$\{(v)|v\in\Omega\}\subseteq A$ a factor set of $B$ in $A$ with
$(v)=1$ if $f_v=y^{\epsilon}y^{-\epsilon}-1, \ y\in Y, \
\epsilon=\pm 1$. Then, for the tower order defined as before, $S$ is
a Gr\"{o}bner-Shirshov bases for $G$ if and only if for any
$v\in\Omega, \ a\in A$ and any composition $(f_{v_1},f_{v_2})_{_w}$
of $R$ in $k\langle Y \rangle$,
\begin{eqnarray}\label{e4}
(v) a^{v}= a^{h_v}(v) \ \ \mbox{ and } \ \ \
\xi_{(v_1,v_2)_{_w}}(v)=\zeta_{(v_1,v_2)_{_w}}(v)
\end{eqnarray}
hold in $A$, where $\xi_{(v_1,v_2)_{_w}}(v)$ and
$\zeta_{(v_1,v_2)_{_w}}(v)$ are defined by (\ref{e2}). Moreover,
if this is the case, $G$  is a Schreier extension of $A$ by $B$ in
a natural way.
\end{theorem}
\noindent {\bf Remark.} \  We call (\ref{e4}) the extension
conditions of $A$ by $B$.

\noindent {\bf Proof.} \  If $S$ is a Gr\"{o}bner-Shirshov bases in
$k\langle A_1\cup Y\rangle$ then by Composition-Diamond Lemma, the
group $\bar A$ is the same as $A$. Then, (\ref{e6}) and (\ref{e3})
imply (\ref{e4}). Conversely, it is trivial that (\ref{e4}) implies
(\ref{e6}) and (\ref{e3}). Now, we need only to prove the last
statement of the theorem that $G$ is a Schreier extension of $A$ by
$B$ in a natural way. By Composition-Diamond Lemma, each element in
$G$ can be uniquely written as $ba$, where $a\in A$ and $b\in Y^*$
is $R$-irreducible. In $G$, the multiplication is defined as
follows: for any $a,a'\in A, \ b,b'\in Y^*$,
$$
ba\cdot b'a'=[bb']\theta_{bb'}(v)a^{b'}a'
$$
where $bb'\equiv [bb'] \ mod(R), \ [bb']$ is $R$-irreducible and
$\theta_{bb'}(v)\in A$ is a function of $\{(v)|v\in\Omega\}$. From
this it follows that $A$ is normal subgroup of $G$ and $G/A\cong B$
as groups. \ \ $\square$

\begin{lemma}\label{sgt2}
Let the notations be as before. Let the group $C$ be a Schreier
extension of $A$ by $B$. Then, there exist a factor set
$\{(v)|v\in\Omega\}$ of $B$ in $A$  with $(v)=1$ if
$f_v=y^{\epsilon}y^{-\epsilon}-1, \ y\in Y, \ \epsilon=\pm 1$ and
$\{a^y|y\in Y, \ A\rightarrow A, \ a\mapsto a^y \mbox{ is an
isomorphism}\}$ such that (\ref{e4}) holds. Moreover, if this is the
case, $C\cong G$, where $G$ is defined in Theorem \ref{sgt1}.
\end{lemma}

\noindent{\bf Proof.} \ Let $\sigma: \ B\rightarrow C/A, \ y\mapsto
g_{y}A$ be an isomorphism. Since, in $B$, $v=h_v$ for any
$v\in\Omega$, we have $g_{v}A=g_{_{h_v}}A$. Then, for any
$v\in\Omega$, there exists a unique $(v)\in A$ such that
$g_{v}=g_{_{h_v}}\cdot(v)$. Thus, we obtain a factor set
$\{(v)|v\in\Omega\}$ of  $B$ in $A$.

For any $y\in Y$, let $\psi_y: \ A\rightarrow A, \ a\mapsto
a^{y}=g_{y}^{-1}ag_{y}$. Then, it is clear that $\psi_y$ is an
automorphism of $A$ because $A\lhd C$.

 Since $\sigma$ is an
isomorphism, each element $c\in C$ can be uniquely expressed as
$c=g_ba$, where $a\in A$ and $b\in Y^*$ is $R$-irreducible. Hence,
by Composition-Diamond Lemma, $S$ is a Gr\"{o}bner-Shirshov bases in
$k\langle A_1\cup Y\rangle$ and so (\ref{e4}) holds.

Let $c'=g_{b'}a'\in C$. Then
$$
cc'=g_bag_{b'}a'=g_bg_{b'}a^{b'}a'=g_{bb'}\theta_{bb'}(v)a^{b'}a'
$$
where $\theta_{bb'}(v)\in A$ is a function of $\{(v)|v\in\Omega\}$.
Now, it is easy to see that $C\rightarrow
G=E(A,Y,a^y,(v))=sgp\langle A_1\cup Y| S\rangle$ by $\tau(g_b)=b, \
\tau(a)=a$ is an isomorphism. \ \ $\square$

\ \

By using Theorem \ref{sgt1} and Lemma \ref{sgt2}, we have the
following theorem which characterizes completely the Schreier
extensions of groups.

\begin{theorem}\label{sgt4}
Let  $A,B$ be groups, $B=sgp\langle Y|R \rangle$, where
$R=\{v=h_v|v\in\Omega\}$ is a  minimal Gr\"{o}bner-Shirshov bases
for $B$ and $v$ the leading term of the polynomial $v-h_v$ in
$k\langle Y \rangle$. Then, a group $C$ is a Schreier extension of
$A$ by $B$ if and only if there exist $\{a^y|y\in Y, \ A\rightarrow
A, \ a\mapsto a^y \mbox{ is an isomorphism}\}$ and a factor set
$\{(v)|v\in\Omega\}$ of  $B$ in $A$  with $(v)=1$ if
$f_v=y^{\epsilon}y^{-\epsilon}-1, \ y\in Y, \ \epsilon=\pm 1$
 such that (\ref{e4}) holds. Moreover, if this is the
case, $C\cong G=E(A,Y,a^y,(v))=sgp\langle A_1\cup Y| S\rangle$,
where $A_1=A\backslash\{1\}, \ S=\{aa'=[aa'], \ v=h_v\cdot(v), \
ay=ya^y |v\in\Omega, \ a,a'\in A_1, \ y\in Y\}$.
\end{theorem}

\ \

\noindent {\bf Remark.} \ In the above theorems, if the group $A$ is
also presented by generators and relations: $A=sgp\langle X|R_A
\rangle$, where $R_A$ is a Gr\"{o}bner-Shirshov bases for $A$ with
the deg-lex order $<_A$ on $X^*$. Then, by replacing $A_1$ with $X$,
$\{aa'=[aa']|a,a'\in A_1\}$ with $R_A$ and $a$ with $x \ (x\in X)$,
the results hold.

\ \

A special case of Theorem \ref{sgt4} is when the group $B$ is
presented as $B=sgp\langle B_1| R\rangle$, where
$B_1=B\backslash\{1\}, \ R=\{bb'=[bb']| b,b'\in B_1\}$. It is
obvious that, for a deg-lex order, $R$ is a  Gr\"{o}bner-Shirshov
bases for $B$. Now, if this is the case, it is easy to check, by
(\ref{e6}) and the algorithm (\ref{e2}), that the condition
(\ref{e4}) is: for any $b,b',b''\in B, a\in A$,
\begin{eqnarray}\label{e7}
(b,b')a^{bb'}=a^{[bb']}(b,b') \ \ \mbox{ and } \ \ \
(b,b'b'')(b',b'')=(bb',b'')(b,b')^{b''}
\end{eqnarray}

Then, we have the following theorem.

\begin{theorem}\label{sgt5}
Let  $A,B$ be groups. Let
$$
G=E(A,B,a^b,(b,b'))=sgp\langle A_1\cup B_1| S\rangle
$$
where $A_1=A\backslash\{1\}, \ B_1=B\backslash\{1\}, \
S=\{aa'=[aa'], \ bb'=[bb'](b,b'), \ ab=ba^b |a,a'\in A_1, \ b,b'\in
B_1\}$, $\{(b,b')|b,b'\in B\}$ a factor set of $B$ in $A$ with
$(b,b')=1$ if $b'=b^{-1}$ and $b: \ A\rightarrow A, \ a\mapsto a^b$
an automorphism. Then the following statements hold.
\begin{enumerate}
\item[1.] \ For the tower order defined as before, $S$ is a Gr\"{o}bner-Shirshov
bases for $G$ if and only if (\ref{e7}) holds. Moreover, if this is
the case, $G$  is a Schreier extension of $A$ by $B$ in a natural
way.
\item[2.] \  A group $C$ is a Schreier extension of $A$ by $B$ if
and only if there exist a factor set $\{(b,b')|b,b'\in B\}$ of  $B$
in $A$ and $\{b: \ A\rightarrow A, \ a\mapsto a^b \mbox{ is an
isomorphism}\}$ such that (\ref{e7}) holds. Moreover, if this is the
case, $C\cong G$.
\end{enumerate}
\end{theorem}

\noindent{\bf Remark}: Theorem \ref{sgt5}, 2 is Schreier's theorem
in  \cite{s1}.

\ \

Now, we give some applications of the above results.

Let $A,B, G$ be groups. The previous theorems give an answer to how
to find the conditions which makes $G$ to be a Schreier extension of
$A$ by $B$. In fact, it is essential to find the extension
conditions (\ref{e4}), that is, to find functions
$\xi_{(v_1,v_2)_{_w}}(v), \ \zeta_{(v_1,v_2)_{_w}}(v)$ which are
given by the algorithm (\ref{e2}).

As results, by using the extension conditions (\ref{e4}), let us
give some examples. We give the characterization of the extension of
$A$ by $B$ when the $B$ is cyclic group and free commutative group,
respectively.

\begin{theorem} \ (\cite{h} Theorem 15.3.4)
Let $A,B, G$ be groups and $B=sgp\langle x|x^n=1 \rangle$ the cyclic
group of the order $n$. Then, $G$ is isomorphic to an extension of
$A$ by $B$ if and only if there exist $a_0\in A$ and an automorphism
$\varphi$ of $A$ such that
$$
a^{\varphi^n}=a_0^{-1}aa_0 \ \  \ and \ \ \ a_0^{\varphi}=a_0
$$
Moreover, if this is the case, $G\cong E(A,x,a^{x},(a_0))=sgp\langle
A_1\cup \{x\}| S\rangle$, where $A_1=A\backslash\{1\}, \
S=\{aa'=[aa'], \ x^n=a_0, \ ax=xa^{\varphi} |a,a'\in A_1\}$.
\end{theorem}
\noindent {\bf Proof.} \ Clearly, $R=\{x^n=1\}$ is  a
Gr\"{o}bner-Shirshov bases for $B$. We need only to consider the
composition in $S$: $(x^n-a_0,x^n-a_0)_w, \ w=x^{n+1}$. Since
$(x^n-a_0,x^n-a_0)_w=-a_0x+xa_0\equiv -xa_0^x+xa_0=x(a_0-a_0^x)$, we
obtain the extension condition mentioned in the theorem, where
$a^x=a^{\varphi}$. Now, by Theorem \ref{sgt4}, the result follows. \
\ $\square$

\begin{theorem}
Let $A,B, G$ be groups, $X=\{x_i|i\in I\}, \ I$ a total ordered set.
We extend the order in $I$ to the set $X\cup X^{-1}$: for any
$i,j\in I$, $x_i^{-1}<x_i<x_j^{-1}<x_j$ \ if \ $i<j$ and then we
order $(X\cup X^{-1})^*$ by the deg-lex order. Let $B=sgp\langle
X\cup X^{-1}|R\rangle$, where $R=\{x_p^{\varepsilon}
x_q^{\delta}=x_q^{\delta}x_p^{\varepsilon}, \
x_q^{\varepsilon}x_q^{-\varepsilon}=1| \varepsilon,\delta=\pm1, \
p>q, \ p,q\in I\}$. Then, $G$ is isomorphic to an extension of $A$
by $B$ if and only if there exist a factor set
$\{(x_p^{\varepsilon}, x_q^{\delta})|\varepsilon,\delta=\pm1, \ p>q,
\ p,q\in I \}$ of $B$ in $A$ and automorphisms $x: \ A\rightarrow A,
\ a\mapsto a^x, \ x\in X\cup X^{-1}$ such that for any $p,q,r\in I,
\ p>q>r, \  \varepsilon,\delta,\gamma=\pm1, \ a\in A$, the following
equations hold.
\begin{enumerate}
\item[1.] \ $(x_p^{\varepsilon}, x_q^{\delta})(x_p^{\varepsilon},
 x_r^{\gamma})^{x_q^{\delta}}(x_q^{\delta},x_r^{\gamma})
=(x_q^{\delta},x_r^{\gamma})^{x_p^{\varepsilon}}(x_p^{\varepsilon},
 x_r^{\gamma})(x_p^{\varepsilon}, x_q^{\delta})^{x_r^{\gamma}}$

\item[2.] \ $(x_p^{\varepsilon},x_q^{-\delta})(x_p^{\varepsilon},x_q^{\delta})^{x_q^{-\delta}}=1$

\item[3.] \ $(x_p^{\varepsilon},x_q^{\delta})^{x_p^{-\varepsilon}}(x_p^{-\varepsilon},x_q^{\delta})=1$

\item[4.] \ $a^{x_q^{\delta}x_p^{\varepsilon}}(x_p^{\varepsilon}, x_q^{\delta})
=(x_p^{\varepsilon}, x_q^{\delta})a^{x_p^{\varepsilon}
x_q^{\delta}}$

\item[5.] \ $a^{x_r^{\varepsilon}x_r^{-\varepsilon}}=a$

\end{enumerate}
Moreover, if this is the case, $G\cong sgp\langle A_1\cup X\cup
X^{-1}| S\rangle$, where $A_1=A\backslash\{1\}, \
 S=\{aa'=[aa'],ax^{\varepsilon}=xa^{x^{\varepsilon}},x_p^{\varepsilon}
x_q^{\delta}=x_q^{\delta}x_p^{\varepsilon}\cdot(x_p^{\varepsilon},
x_q^{\delta}), \ x_q^{\varepsilon}x_q^{-\varepsilon}=1|
\varepsilon,\delta=\pm1, \ p>q, \ p,q\in I, \ x\in X, a,a'\in
A_1\}$.
\end{theorem}
\noindent {\bf Proof.} \ It is easy to check that for the given
order, $R$ is a Gr\"{o}bner-Shirshov bases for $B$. Then, $B$ is the
free abelian group generated by $X$. The possible compositions in
$R$ are related to the following ambiguities:
$w_1=x_p^{\varepsilon}x_q^{\delta}x_r^{\gamma}, \
w_2=x_p^{\varepsilon}x_q^{\delta}x_q^{-\delta}, \
w_3=x_p^{\varepsilon}x_p^{-\varepsilon}x_q^{\delta}, \
w_4=x_r^{\varepsilon}x_r^{-\varepsilon}x_r^{\varepsilon}, p,q,r\in
I, \ p>q>r$. Then, in $S$, by calculating the corresponding
compositions and by noting that the factors
$(x_q^{\varepsilon},x_q^{-\varepsilon})=1, \ q\in I$, we obtain the
equations 1-3, respectively. The equations 4 and 5 follow from the
formula (\ref{e6}). Now, by Theorem \ref{sgt4}, the result follows.
\ \ $\square$

\section{Structure of Schreier extension  by an HNN group}

In this section, we give the structure of Schreier extension of
group $A$ by $B$, where the group $B$ is an HNN-extension of a
group.

A Gr\"{o}bner-Shirshov bases for HNN-extension of group is given in
\cite{cz} by using the Generalized Composition-Diamond Lemma (Lemma
\ref{l3.1}) with a special total order (not monomial). Thus, we can
not use directly the results of the above section to find the
extension conditions for HNN-extension of group.

\begin{definition} (\cite{hnn},\cite{n52},\cite{n54},\cite{n55}) \
Let $H$ be a group and let $C$ and $D$ be subgroups of $H$ with
$\phi:C\rightarrow D$ an isomorphism. Then the HNN-extension of
$H$ relative to $C\ ,D$ and $\phi$ is the group
$$
B=gp\langle H,t;t^{-1}ct=\phi(c),c\in C\rangle.
$$
\end{definition}

Let
$$
H/C=\{g_i C|i\in I\}, \ H/D=\{h_j D|j\in J\},
$$
where $\{g_i|i\in I\}$ and $\{h_j|j\in J\}$ are the coset
representatives of $C$ and $D$ in $H$, respectively. Then, for any
$h\in H$, there exist uniquely $h_{_C}\in \{g_i|i\in I\}, h_{_D}\in
\{h_j|j\in J\}, c_{_h}\in C,d_{_h}\in D$ such that
$h=h_{_C}c_{_h}=h_{_D}d_{_h}$.

In \cite{cz}, the group $B$ is presented as
$$
B=sgp\langle H_1,t,t^{-1}|R\rangle,
$$
where $H_1=H\backslash\{1\}, \ R=\{h h^{'}=[h h^{'}],h t=h_{_C}
t\phi(c_{_h}),ht^{-1}=h_{_D} t^{-1}\phi^{-1}(d_{_h}),
t^{\varepsilon}t^{-\varepsilon}=1| h,h'\in H_1, \
\varepsilon=\pm1\}$.

Let $H_1=\{h_{\alpha}|\alpha\in\Omega\}$ and $\Omega, I,J$ the
linearly ordered sets.

Let $A$ be a group,
$$
G=sgp\langle A_1\cup H_1\cup\{t,t^{-1}\}| S\rangle
$$
where $A_1=A\backslash\{1\}$ and $S$ consists of the following
relations:
\begin{enumerate}
\item[(3.1)] $aa'=[aa']$
\item[(3.2)] $ay=ya^y$
\item[(3.3)] $hh'=[hh'](h,h')$
\item[(3.4)] $h t=h_{_C} t\phi(c_{_h})(h,t)$
\item[(3.5)] $ht^{-1}=h_{_D} t^{-1}\phi^{-1}(d_{_h})(h,t^{-1})$
\item[(3.6)] $t^{\varepsilon}t^{-\varepsilon}=1$
\end{enumerate}
where $a,a'\in A_1, \ y\in H_1\cup\{t,t^{-1}\}, \ h,h'\in H_1, \
\varepsilon=\pm1, \ \{(h,h'),(h,t^{\varepsilon})|h,h'\in H_1, \
\varepsilon=\pm1\}\subseteq A$ a factor set of $B$ in $A$ with
$(h,h')=1$ if $h'=h^{-1}$ and $\psi_y: \ A\rightarrow A, \ a\mapsto
a^y$ an automorphism.

The following lemma comes from \cite{cz}.

\begin{lemma}\label{l3.1}
\ (\cite{cz}, Lemma 2.1) \ (Generalized Composition-Diamond Lemma) \
Let $S\subseteq k\langle X\rangle, \ A=k\langle X|S\rangle$ and
$``<"$ a total order on $X^*$ such that
\begin{enumerate}
\item[(I)] \ $\overline{asb}=a\bar{s}b$ for any $a,b\in X^*, \ s\in S$;
\item[(II)] \ for each composition $(s_1,s_2)_w$ in $S$, there exists a presentation
$$
(s_1,s_2)_w=\sum_{i}\alpha_{i}a_it_ib_i, \ a_i\bar{t_i}b_i<w, \ \
\mbox{ where } \ t_i\in S, \ a_i,b_i\in X^*, \ \alpha_{i}\in k
$$
 such that for any $c,d\in X^*$, we have
\begin{eqnarray*}
ca_i\bar{t_i}b_id<cwd
\end{eqnarray*}
\end{enumerate}
Then, the following statements hold.
\begin{enumerate}
\item[(i)] \ $S$ is a Gr\"{o}bner-Shirshov basis.
\item[(ii)] \ For any $f\in k\langle X\rangle, \
0\neq f\in Id(S)\Rightarrow \bar{f}=a\bar{s}b$ for some $s\in S, \
a,b\in X^*$.
\item[(iii)] \
The set
$$
Irr(S)=\{u\in X^*|u\neq a\bar{s}b,s\in S,a,b\in X^*\}
$$
is a basis of the algebra $A$.
\end{enumerate}
\end{lemma}

In \cite{cz}, a special total order (not monomial) $<_B$ on
$(H_1\cup\{t,t^{-1}\})^*$ is defined  such that $R$ satisfies the
conditions (I) and (II) in the Lemma \ref{l3.1}, and then, with the
order $<_B$, $R$ is a Gr\"{o}bner-Shirshov bases for $B$. Now, we
extend $<_B$ on $(H_1\cup\{t,t^{-1}\})^*$ in \cite{cz}
 to one on $(A_1\cup H_1\cup\{t,t^{-1}\})^*$ such that, with
such an order, $S$ satisfies the conditions (I) and (II) in the
Lemma \ref{l3.1} under some conditions.

Let $A_1$ be a total ordered set. We order the set $(A_1\cup
H_1\cup\{t,t^{-1}\})^*$ by the following steps (also see \cite{cz}).

Step 1:  order the set $H$ in three different ways:
\begin{enumerate}
\item[$(i)_1$] \ Let $1<h_{_\alpha}<h_{_\beta}<\cdots \ (\alpha<\beta)$
be ordered by $\Omega$. Then we denote this order by $(H,>)$ and
call it an absolute order.
\item[$(ii)_1$] \ For any $h,h^{'}\in H$, suppose that $h=h_{_C}
c_h,h^{'}=h^{'}_{_C} c_{h^{'}}$. Then $h>_{_C}h^{'}$ if and only if
$(h_{_C},c_h)>(h^{'}_{_C},c_h^{'})$ which is ordered
lexicographically: elements $h_{_C},h^{'}_{_C}$ by $I$, elements
$c_h,c_{h^{'}}$ by $(i)_1$. We denote this order by $(H,>_{_C})$ and
call it the $C$-order.
\item[$(iii)_1$] \ For any $h,h^{'}\in H$, suppose that
$h=h_{_D} d_h,h^{'}=h'_{_D} d_{h^{'}}$. Then $h>_{_D}h^{'}$ if and
only if $(h_{_D},d_h)>(h'_{_D},d_{h^{'}})$ which is ordered
lexicographically: elements $h_{_D},h^{'}_{_D}$ by $J$, elements
$d_h,d_{h^{'}}$ by $(i)_1$. We denote this order by $(H,>_{_D})$ and
call it the $D$-order.
\end{enumerate}

Step 2: order the set $(A_1\cup H_1)^{*}$ in three different ways
too.

Each element in $(A_1\cup H_1)^{*}$ has a unique form $e=e_1
h_1\cdots e_nh_n e_{n+1}$, where each ${e_i}\in A_1^*, \ h_i\in H_1,
\ n\geq 0$. Suppose that $e'=e_1' h_1'\cdots e_l'h_l'e_{l+1}'\in
(A_1\cup H_1)^{*}$. Let
\begin{eqnarray*}
&&wt_h(e)=(n,h_1,\cdots,h_n,e_1,\cdots,e_n,e_{n+1})\\
&&wt_h(e')=(l,h_1',\cdots,h_l',e_1',\cdots,e_l',e_{l+1}')
\end{eqnarray*}
We define $e<e' \ (e<_{_C}e', \ e<_{_D}e')$ if $wt_h(e)<wt_h(e') \
(wt_h(e)<_{_C}wt_h(e'), \ wt_h(e)<_{_D}wt_h(e'))$ lexicographically,
using the order of natural numbers and the following orders:
\begin{enumerate}
\item[$(i)_2$]\ The absolute order $((A_1\cup H_1)^{*},\leq)$: to
compare $h_i,h_i'$ by $(i)_1$ and $e_i,e_i'$ by the deg-lex order on
$A_1$.

\item[$(ii)_2$]\ The $C$-order $((A_1\cup H_1)^{*},\leq_{_C})$:
firstly to compare  $h_1,\cdots,h_{n-1}, \ (n\geq1)$,
lexicographically by absolute order $(i)_1$,  secondly for the last
element $h_n$ by $C$-order $(ii)_1$ and finally for $e_i,e_i'$ by
the deg-lex order on $A_1$.
\item[$(iii)_2$]\ The $D$-order $((A_1\cup H_1)^{*},\leq {_D})$ is similar to
$(ii)_2$ by replacing $>_{_C} \ \mbox{with} \ >_{_D}$.
\end{enumerate}

Step 3: order the set $(A_1\cup H_1\cup\{t,t^{-1}\})^{*}$.

Each element in $(A_1\cup H_1\cup\{t,t^{-1}\})^{*}$ has a unique
form $u=u_1 t^{{\varepsilon}_1}\cdots u_kt^{{\varepsilon}_k}
u_{k+1}$, where each ${u_i}\in (A_1\cup H_1)^{*},
{\varepsilon_i}=\pm1, \ k\geq 0$. Suppose that $v=v_1
t^{{\delta}_1}\cdots v_lt^{{\delta}_l}v_{l+1}\in(A_1\cup
H_1\cup\{t,t^{-1}\})^{*}$. Let
\begin{eqnarray*}
&&wt(u)=(k,t^{\varepsilon_1},\cdots,t^{\varepsilon_k},wt_h(u_1),
\cdots,wt_h(u_k),wt_h(u_{k+1}))\\
&&wt(v)=(l,t^{\delta_1},\cdots,t^{\delta_l},wt_h(v_1),\cdots,wt_h(v_l),wt_h(v_{l+1}))
\end{eqnarray*}
 We define $u\succ v$ if $wt(u)>wt(v)$ lexicographically, using the
order of natural numbers and the following orders:
\begin{enumerate}
\item[(a)]\ $t>t^{-1}$
\item[(b)]\ $wt_h(u_i)>_{_C}wt_h(v_i) \ \mbox{if} \ \varepsilon_i=1, \ 1\leq i\leq k$
\item[(c)]\ $wt_h(u_i)>_{_D}wt_h(v_i) \ \mbox{if} \ \varepsilon_i=-1, \ 1\leq i\leq k$
\item[(d)]\ $wt_h(u_{k+1})> wt_h(v_{l+1}) \ (k=l)$, the absolute order $(i)_2$
\end{enumerate}

In $S$, the possible compositions  are related to ambiguities:
\begin{eqnarray*}
&&w_1=hh'h'', \ w_2=hh't, \ w_3=hh't^{-1}, \ w_4=htt^{-1}, \
w_5=ht^{-1}t\\
&&w_6=ahh', \ w_7=aht, \ w_8=aht^{-1}, \
w_{9}=at^{\varepsilon}t^{-\varepsilon}, \
w_{10}=t^{\varepsilon}t^{-\varepsilon}t^{\varepsilon}
\\
&&w_{11}=aa'a'', \ w_{12}=aa'y
\end{eqnarray*}
where $h,h',h''\in H_1, \ a,a',a''\in A_1, \ y\in
H_1\cup\{t,t^{-1}\}, \ {\varepsilon_i}=\pm1$. By using algorithm
(\ref{e4}), we have the following conditions ($w_i\leftrightarrow
(hi), \ 1\leq i\leq 9$) in $A$:
\begin{enumerate}
\item[(h1).] \ $(h,h'h'')(h',h'')=(hh',h'')(h,h')^{h''}$

\item[(h2).] \ $(hh',t)(h,h')^t=(hh'_{_C},t)^{\phi(c_{_{h'}})}(h',t)$

\item[(h3).] \ $(hh',t^{-1})(h,h')^{t^{-1}}=
(hh'_{_D},t^{-1})^{\phi^{-1}(d_{_{h'}})}(h',t^{-1})$

\item[(h4).] \ $(\phi(c_{_{h}}) ,t^{-1})(h,t)^{t^{-1}}=1$

\item[(h5).] \ $(\phi^{-1}(d_{_{h}}) ,t)(h,t^{-1})^{t}=1$

\item[(h6).] \ $a^{hh'}(h,h')=(h,h')a^{[hh']}$

\item[(h7).] \ $a^{^{h_{_C}t\phi(c_{_h})}}(h,t)=(h,t)a^{ht}$

\item[(h8).] \ $a^{^{h_{_D} t^{-1}\phi^{-1}(d_{_h})}}(h,t^{-1})=(h,t^{-1})a^{ht^{-1}}$

\item[(h9).] \ $a^{t^{\varepsilon}t^{-\varepsilon}}=a$

\end{enumerate}

\begin{lemma}\label{l3.2}
Let the notations as above and the conditions (I) and (II) be in
Lemma \ref{l3.1}. Then, for $S$ and the order $\succ$, in $k\langle
A_1\cup H_1\cup\{t,t^{-1}\}\rangle$, the following statements hold.
\begin{enumerate}
\item[(i)] \ (I) holds. \item[(ii)] \  $S$ and $\succ$ satisfy
(II) if and only if
 $(hi), \ 1\leq i\leq 9$ hold.
\end{enumerate}
\end{lemma}
\noindent{\bf Proof.} \ Let $c,d\in (A_1\cup H_1\cup\{t,t^{-1}\})^*,
\ c=c_1 t^{{\varepsilon}_1} \cdots c_kt^{{\varepsilon}_k} c_{k+1}, \
d=d_1 t^{{\delta}_1}\cdots d_lt^{{\delta}_l}d_{l+1}$ where each
$c_i,d_i\in (A_1\cup H_1)^{*}, \varepsilon_i,\delta_i=\pm1, \
k,l\geq 0$. Let $c_{k+1}=a_1h_1\cdots a_ph_pa_{p+1}, \
d_1=a_1'h_1'\cdots a_q'h_q'a_{q+1}', \ a_i,a_i'\in A_1^*, \
h_i,h_i'\in H_1,p,q\geq 0$.

(i) To prove $\overline{csd}=c\bar{s}d$ for any $s\in S$, it
suffices to prove:
\begin{enumerate}
\item[(I1)] $caa'd\succ c[aa']d$
\item[(I2)] $cayd\succ cya^yd$
\item[(I3)] $chh'd\succ c[hh'](h,h')d$
\item[(I4)] $ch td\succ ch_{_C} t\phi(c_{_h})(h,t)d$
\item[(I5)] $cht^{-1}d\succ ch_{_D} t^{-1}\phi^{-1}(d_{_h})(h,t^{-1})d$
\item[(I6)] $ct^{\varepsilon}t^{-\varepsilon}d\succ cd$
\end{enumerate}

We only check (I2) and (I4). For (I1), (I3), (I5) and (I6), the
proof are similar.

For (I2), we have three cases to consider: $y=t, \ y=t^{-1}, y=h\in
H_1$. Suppose that $y=h$. Then
\begin{eqnarray*}
wt(cahd)&=&(k+l,t^{\varepsilon_1},\cdots,t^{\varepsilon_k},
t^{\delta_1},\cdots,t^{\delta_l},wt_h(c_1),\cdots,wt_h(c_k),\\
&&wt_h(c_{k+1}ahd_1), wt_h(d_2),\cdots,wt_h(d_{l+1}))\\
wt(cha^hd)&=&(k+l,t^{\varepsilon_1},\cdots,t^{\varepsilon_k},
t^{\delta_1},\cdots,t^{\delta_l},wt_h(c_1),\cdots,wt_h(c_k),\\
&&wt_h(c_{k+1}ha^hd_1), wt_h(d_2),\cdots,wt_h(d_{l+1}))
\end{eqnarray*}
Also,
\begin{eqnarray*}
wt_h(c_{k+1}ahd_1)&=&(p+q+1,h_1,\cdots,h_p,h,h_1',\cdots,h_q',a_1,\cdots,\\
&&a_p,a_{p+1}a,a_1',\cdots)\\
wt_h(c_{k+1}ha^hd_1)&=&(p+q+1,h_1,\cdots,h_p,h,h_1',\cdots,h_q',a_1,\cdots,\\
&&a_p,a_{p+1},a^ha_1',\cdots)
\end{eqnarray*}
Since $a_{p+1}a>a_{p+1}$, we have
$wt_h(c_{k+1}ahd_1)>wt_h(c_{k+1}ha^hd_1)$ and hence, $wt(cahd)\succ
wt(cha^hd)$.

Suppose that $y=t$. Then
\begin{eqnarray*}
wt(catd)&=&(k+l+1,t^{\varepsilon_1},\cdots,t^{\varepsilon_k},t,
t^{\delta_1},\cdots,t^{\delta_l},wt_h(c_1),\cdots,wt_h(c_k),\\
&&wt_h(c_{k+1}a), wt_h(d_1),\cdots, wt_h(d_l),wt_h(d_{l+1}))\\
wt(cta^td)&=&(k+l+1,t^{\varepsilon_1},\cdots,t^{\varepsilon_k},t,
t^{\delta_1},\cdots,t^{\delta_l},wt_h(c_1),\cdots,wt_h(c_k),\\
&&wt_h(c_{k+1}), wt_h(a^td_1),\cdots, wt_h(d_l),wt_h(d_{l+1}))
\end{eqnarray*}
Also,
\begin{eqnarray*}
wt_h(c_{k+1}a)&=&(p,h_1,\cdots,h_p,a_1,\cdots,a_p,a_{p+1}a)\\
wt_h(c_{k+1})&=&(p,h_1,\cdots,h_p,a_1,\cdots,a_p,a_{p+1})
\end{eqnarray*}
Since $a_{p+1}a>a_{p+1}$, we have $wt_h(c_{k+1}a)>wt_h(c_{k+1})$ and
hence, $wt(catd)\succ wt(cta^td)$.

Similarly, for $y=t^{-1}$, $wt(cayd)\succ wt(cya^yd)$ holds.

For (I4), we have
\begin{eqnarray*}
wt(chtd)&=&(k+l+1,t^{\varepsilon_1},\cdots,t^{\varepsilon_k},t,
t^{\delta_1},\cdots,t^{\delta_l},wt_h(c_1),\cdots,wt_h(c_k),\\
&&wt_h(c_{k+1}h), wt_h(d_1),\cdots, wt_h(d_l),wt_h(d_{l+1}))\\
wt(ch_{_C}
t\phi(c_{_h})(h,t)d)&=&(k+l+1,t^{\varepsilon_1},\cdots,t^{\varepsilon_k},t,
t^{\delta_1},\cdots,t^{\delta_l},wt_h(c_1),\cdots,wt_h(c_k),\\
&&wt_h(c_{k+1}h_{_C}), wt_h(\phi(c_{_h})(h,t)d_1),\cdots,
wt_h(d_l),wt_h(d_{l+1}))
\end{eqnarray*}
Also,
\begin{eqnarray*}
wt_h(c_{k+1}h)&=&(p+1,h_1,\cdots,h_p,h,a_1,\cdots,a_p,a_{p+1},1)\\
wt_h(c_{k+1}h_{_C})&=&(p+1,h_1,\cdots,h_p,h_{_C},a_1,\cdots,a_p,a_{p+1},1)
\end{eqnarray*}
Since $\varepsilon=1$ and $h_{_C}<_{_C}h$, we have
$wt_h(c_{k+1}h_{_C})<_{_C}wt_h(c_{k+1}h)$ and hence, $chtd\succ
ch_{_C} t\phi(c_{_h})(h,t)d$.

Therefore, we show (i).

(ii) It is easy to check that, for
$w_{10}=t^{\varepsilon}t^{-\varepsilon}t^{\varepsilon}, \
w_{11}=aa'a''$ and $w_{12}=aa'y$, the corresponding composition in
$S$ satisfies (II).

Now, we prove that for  $w_i, \ 1\leq i\leq 9$, the corresponding
composition in $S$ satisfies (II) if and only if $(hi)$ holds.

 For  $w_i, \ 1\leq i\leq 9$, the corresponding composition in $S$ are,
respectively,
\begin{enumerate}
\item[(II1)] $(hh'-[hh'](h,h'),h'h''-[h'h''](h',h''))_{w_1}, \ w_1=hh'h''$
\item[(II2)] $(hh'-[hh'](h,h'),h't-h'_{_C} t\phi(c_{_{h'}})(h',t))_{w_2}, \ w_2=hh't$
\item[(II3)] $(hh'-[hh'](h,h'),h't^{-1}-h'_{_D} t^{-1}\phi^{-1}(d_{_{h'}})
(h',t^{-1}))_{w_3}, \ w_3=hh't^{-1}$
\item[(II4)]  $(ht-h_{_C} t\phi(c_{_h})(h,t),
tt^{-1}-1)_{w_4}, \ w_4=htt^{-1}$
\item[(II5)] $(ht^{-1}-h_{_D} t^{-1}\phi^{-1}(d_{_h})(h,t^{-1}),
t^{-1}t-1)_{w_5}, \ w_5=ht^{-1}t$
\item[(II6)] $(ah-ha^h,hh'-[hh'](h,h'))_{w_6}, \ w_6=ahh'$
\item[(II7)] $(ah-ha^h,ht-h_{_C}t\phi(c_{_h})(h,t))_{w_7}, \ w_7=aht$
\item[(II8)] $(ah-ha^h,ht^{-1}-h_{_D}
t^{-1}\phi^{-1}(d_{_h})(h,t^{-1}))_{w_8}, \ w_8=aht^{-1}$
\item[(II9)] $(at^{\varepsilon}-t^{\varepsilon}a^{t^{\varepsilon}},
t^{\varepsilon}t^{-\varepsilon}-1)_{w_{9}}, \
w_{9}=at^{\varepsilon}t^{-\varepsilon}$
\end{enumerate}

We will check (II2), (II4) and (II7).

 For (II2), $ w_2=hh't$, by
noting that $[hh'_{_C}]_{_C}=[hh']_{_C}$ and
$c_{_{[hh'_{_C}]}}c_{_{h'}}=c_{_{[hh']}}$, we have
\begin{eqnarray*}
&&(hh'-[hh'](h,h'),h't-h'_{_C} t\phi(c_{_{h'}})(h',t))_{w_2}\\
&=&-[hh'](h,h')t+hh'_{_C}t\phi(c_{_{h'}})(h',t)\\
&=&-[hh']((h,h')t-t(h,h')^t)-([hh']t-[hh']_{_C}t\phi(c_{_{[hh']}})(hh',t))(h,h')^t\\
&&+(hh'_{_C}t-[hh'_{_C}]_{_C}t\phi(c_{_{[hh'_{_C}]}})(hh'_{_C},t))\phi(c_{_{h'}})(h',t)\\
&&+[hh'_{_C}]_{_C}t\phi(c_{_{[hh'_{_C}]}})((hh'_{_C},t)\phi(c_{_{h'}})-
\phi(c_{_{h'}})(hh'_{_C},t)^{\phi(c_{_{h'}})})(h',t)\\
&&+[hh'_{_C}]_{_C}t\phi(c_{_{[hh'_{_C}]}})\phi(c_{_{h'}})
(hh'_{_C},t)^{\phi(c_{_{h'}})}(h',t)-
[hh']_{_C}t\phi(c_{_{[hh']}})(hh',t)(h,h')^t\\
&=&-[hh']t_1-t_2(h,h')^t +t_3\phi(c_{_{h'}})(h',t)
+[hh'_{_C}]_{_C}t\phi(c_{_{[hh'_{_C}]}})t_4(h',t)\\
&&+
[hh']_{_C}t\phi(c_{_{[hh']}})((hh'_{_C},t)^{\phi(c_{_{h'}})}(h',t)-(hh',t)(h,h')^t)
\end{eqnarray*}
where
\begin{eqnarray*}
&&t_1=(h,h')t-t(h,h')^t, \ \ t_2=
[hh']t-[hh']_{_C}t\phi(c_{_{[hh']}})(hh',t)\\
&&t_3=hh'_{_C}t-[hh'_{_C}]_{_C}t\phi(c_{_{[hh'_{_C}]}})(hh'_{_C},t),
\ \ t_4=(hh'_{_C},t)\phi(c_{_{h'}})-
\phi(c_{_{h'}})(hh'_{_C},t)^{\phi(c_{_{h'}})}
\end{eqnarray*}
Clearly, $t_1,t_2,t_3,t_4\in S$ and since
\begin{eqnarray*}
wt_h(c_{k+1}hh')&_{_C}\!\!>&wt_h(c_{k+1}[hh'](h,h')), \ wt_h(c_{k+1}[hh']),\\
&&wt_h(c_{k+1}hh'_{_C}), \ wt_h(c_{k+1}[hh'_{_C}]_{_C})
\end{eqnarray*}
we have
\begin{eqnarray*}
&&c[hh']\bar{t_1}d=c[hh'](h,h')td\prec chh'td\\
&&c\bar{t_2}(h,h')^td=c[hh']t(h,h')^td\prec chh'td\\
&&c\bar{t_3}\phi(c_{_{h'}})(h',t)d=chh'_{_C}t\phi(c_{_{h'}})(h',t)d\prec chh'td\\
&&c[hh'_{_C}]_{_C}t\phi(c_{_{[hh'_{_C}]}})\bar{t_4}(h',t)d
=c[hh'_{_C}]_{_C}t\phi(c_{_{[hh'_{_C}]}})(hh'_{_C},t)\phi(c_{_{h'}})(h',t)d\prec
chh'td
\end{eqnarray*}

For (II4),  $w_4=htt^{-1}$, we have
\begin{eqnarray*}
&&(ht-h_{_C} t\phi(c_{_h})(h,t), tt^{-1}-1)_{w_4}\\
&=&-h_{_C} t\phi(c_{_h})(h,t)t^{-1}+h\\
&=&-h_{_C} t\phi(c_{_h})((h,t)t^{-1}-t^{-1}(h,t)^{t^{-1}})\\
&&-h_{_C}t(\phi(c_{_h})t^{-1}-t^{-1}c_{_h}(\phi(c_{_h}),t^{-1}))
(h,t)^{t^{-1}}\\
&&-h_{_C}(tt^{-1}-1)c_{_h}(\phi(c_{_h}),t^{-1})
(h,t)^{t^{-1}}\\
&&-h(\phi(c_{_h}),t^{-1})
(h,t)^{t^{-1}}+h\\
&=&-h_{_C} t\phi(c_{_h})t_1 -h_{_C}tt_2 (h,t)^{t^{-1}}
-h_{_C}t_3c_{_h}(\phi(c_{_h}),t^{-1})(h,t)^{t^{-1}}
\\
&&-h((\phi(c_{_h}),t^{-1}) (h,t)^{t^{-1}}-1)
\end{eqnarray*}
where $t_1=(h,t)t^{-1}-t^{-1}(h,t)^{t^{-1}}, \
t_2=\phi(c_{_h})t^{-1}-t^{-1}c_{_h}(\phi(c_{_h}),t^{-1}), \
t_3=tt^{-1}-1$. Clearly, $t_1,t_2,t_3\in S$ and since
\begin{eqnarray*}
&&wt_h(c_{k+1}h) \ _{_C}\!\!> \ wt_h(c_{k+1}h_{_C}), \ \ wt(chtt^{-1}d)=(k+l+2,\cdots)\\
&& and \  \ wt(ch_{_C}c_{_h}(\phi(c_{_h}),t^{-1})
 (h,t)^{t^{-1}}d)=(k+l,\cdots)
\end{eqnarray*}
we have
\begin{eqnarray*}
&&ch_{_C} t\phi(c_{_h})\bar{t_1}d=ch_{_C} t\phi(c_{_h})(h,t)t^{-1}d\prec chtt^{-1}d\\
&&ch_{_C}t\bar{t_2}(h,t)^{t^{-1}}d=ch_{_C}t\phi(c_{_h})t^{-1}(h,t)^{t^{-1}}d\prec chtt^{-1}d\\
&&ch_{_C}\bar{t_3}c_{_h}(\phi(c_{_h}),t^{-1})(h,t)^{t^{-1}}d=
ch_{_C}tt^{-1}c_{_h}(\phi(c_{_h}),t^{-1})(h,t)^{t^{-1}}d\prec
chtt^{-1}d
\end{eqnarray*}

For (II7), $w_7=aht$, we have
\begin{eqnarray*}
&&(ah-ha^h,ht-h_{_C}t\phi(c_{_h})(h,t))_{w_7}=-ha^ht+ah_{_C}t\phi(c_{_h})(h,t)\\
&=&-h(a^ht-ta^{ht})-(ht-h_{_C}t\phi(c_{_h})(h,t))a^{ht}+
(ah_{_C}t\phi(c_{_h})-h_{_C}t\phi(c_{_h})a^{^{h_{_C}t\phi(c_{_h})}})(h,t)\\
&&+h_{_C}t\phi(c_{_h})(a^{h_{_C}t\phi(c_{_h})}(h,t)-(h,t)a^{ht})\\
&=&-ht_1-t_2a^{ht}+
t_3(h,t)+h_{_C}t\phi(c_{_h})(a^{^{h_{_C}t\phi(c_{_h})}}(h,t)-(h,t)a^{ht})
\end{eqnarray*}
where $t_1=a^ht-ta^{ht}, \ t_2=ht-h_{_C}t\phi(c_{_h})(h,t), \
t_3=ah_{_C}t\phi(c_{_h})-h_{_C}t\phi(c_{_h})a^{^{h_{_C}t\phi(c_{_h})}}\in
S$. Also,
\begin{eqnarray*}
wt(cw_7d)=wt(cahtd)&=&(k+l+1,t^{\varepsilon_1},\cdots,t^{\varepsilon_k},t,
t^{\delta_1},\cdots,t^{\delta_l},wt_h(c_1),\cdots,\\
&&wt_h(c_k),wt_h(c_{k+1}ah),wt_h(d_1),\cdots)\\
wt(ch\bar{t_1}d)=wt(cha^htd)&=&(k+l+1,t^{\varepsilon_1},\cdots,t^{\varepsilon_k},t,
t^{\delta_1},\cdots,t^{\delta_l},wt_h(c_1),\cdots,\\
&&wt_h(c_k),wt_h(c_{k+1}ha^h),wt_h(d_1),\cdots)\\
wt(c\bar{t_2}a^{ht}d)=wt(chta^{ht}d)&=&(k+l+1,t^{\varepsilon_1},\cdots,t^{\varepsilon_k},t,
t^{\delta_1},\cdots,t^{\delta_l},wt_h(c_1),\cdots,\\
&&wt_h(c_k),wt_h(c_{k+1}h),wt_h(a^{ht}d_1),\cdots)\\
wt(c\bar{t_3}(h,t)d)=wt(cah_{_C}t\phi(c_{_h})(h,t)d)
&=&(k+l+1,t^{\varepsilon_1},\cdots,t^{\varepsilon_k},t,
t^{\delta_1},\cdots,t^{\delta_l},wt_h(c_1),\cdots,\\
&&wt_h(c_k),wt_h(c_{k+1}ah_{_C}),wt_h(\phi(c_{_h})(h,t)d_1),\cdots)
\end{eqnarray*}
and
\begin{eqnarray*}
wt_h(c_{k+1}ah)&=&(p+1,h_1,\cdots,h_p,h,a_1,\cdots,a_p,a_{p+1}a,1)\\
wt_h(c_{k+1}ha^h)&=&(p+1,h_1,\cdots,h_p,h,a_1,\cdots,a_p,a_{p+1},a^h)\\
wt_h(c_{k+1}h)&=&(p+1,h_1,\cdots,h_p,h,a_1,\cdots,a_p,a_{p+1},1)\\
wt_h(c_{k+1}ah_{_C})&=&(p+1,h_1,\cdots,h_p,h_{_C},a_1,\cdots,a_p,a_{p+1}a,1)
\end{eqnarray*}
Since $a_{p+1}a>a_{p+1}$ and $h \ _{_C}\!\!>h_{_C}$, we have
$$
wt_h(c_{k+1}ah) \ _{_C}\!\!>wt_h(c_{k+1}ha^h),
wt_h(c_{k+1}h),wt_h(c_{k+1}ah_{_C})
$$
and so
$$
cw_7d\succ ch\bar{t_1}d,c\bar{t_2}a^{ht}d,c\bar{t_3}(h,t)d
$$

All other cases are treated the same as before.

Now, suppose that $(h1)-(h9)$ hold. Then, it is clear that $S$ and
$\succ$ satisfy (II). Conversely, assume that $S$ and $\succ$
satisfy (II). Then, from the previous relations, it follows that
$(h1)-(h9)$ hold in the group $\bar A$, the natural homomorphic
image of A in G. By Generalized Composition-Diamond Lemma, we know
that groups $A$ and $\bar A$ are the same. So, we show (ii).

The proof of the lemma is completed. \ \ $\square$

\ \

By using Lemma \ref{l3.2}, Lemma \ref{l3.1} and the proof of the
theorems in Section 2, we have the following theorem.

\begin{theorem}
Let $A,H, G$ be groups and $B=sgp\langle H_1,t,t^{-1}|R\rangle$
the HNN-extension of the group $H$, where $R=\{h h^{'}=[h h^{'}],h
t=h_{_C} t\phi(c_{_h}),ht^{-1}=h_{_D} t^{-1}\phi^{-1}(d_{_h}),
t^{\varepsilon}t^{-\varepsilon}=1| h,h'\in H_1, \
\varepsilon=\pm1\}, \ H_1=H\backslash\{1\}, \ C,D$ the subgroups
of $H$, $\phi:C\rightarrow D$ an isomorphism,
$h=h_{_C}c_{_h}=h_{_D}d_{_h}, \ h_{_C}\in \{g_i|i\in I\},
h_{_D}\in \{h_j|j\in J\}, c_{_h}\in C,d_{_h}\in D$, $\{g_i|i\in
I\}$ and $\{h_j|j\in J\}$ the coset representatives of $C$ and $D$
in $H$, respectively. Then, $G$ is isomorphic to a Schreier
extension of $A$ by $B$ if and only if there exist a factor set
$\{(h,h'),(h,t^{\varepsilon})|h,h'\in H_1, \
\varepsilon=\pm1\}\subseteq A$ of $B$ in $A$ with $(h,h')=1$ if
$h'=h^{-1}$ and $\psi_b: \ A\rightarrow A, \ a\mapsto a^b \ \
(b\in H_1\cup\{t,t^{-1}\})$ an automorphism such that for any
$h,h',h''\in H, \ a\in A, \ \varepsilon=\pm1$, the equations
(h1)-(h9) hold.

Moreover, if this is the case, $ G\cong sgp\langle A_1\cup
H_1\cup\{t,t^{-1}\}| S\rangle $, where $A_1=A\backslash\{1\}, \
H_1=H\backslash\{1\}, \ S=\{aa'=[aa'],ab=ba^b,h h^{'}=[h
h^{'}](h,h^{'}),ht=h_{_C} t\phi(c_{_h})(h,t),ht^{-1}= h_{_D}
t^{-1}\phi^{-1}(d_{_h})(h,t^{-1}),
t^{\varepsilon}t^{-\varepsilon}=1| a,a'\in A_1,b\in
H_1\cup\{t,t^{-1}\}, h,h'\in H_1, \varepsilon=\pm1\}$.

\end{theorem}

 \ \

\noindent{\bf Acknowledgement}: The author would like to express his
deepest gratitude to Professor L. A. Bokut for his kind guidance,
useful discussions and enthusiastic encouragement  when the author
was visiting Sobolev Institute of Mathematics.

 \ \

\end{document}